\documentclass[10pt]{article}
\usepackage{amssymb,amsmath,amsthm,epsfig}
\usepackage[dvips]{color}
\newtheorem{theorem}{Theorem}

\newtheorem{cor}{Corollary}

\newtheorem{lemma}{Lemma}
\allowdisplaybreaks
\begin{document}

\title{Discrete Green's functions for products of regular graphs}
\author{Robert B.\ Ellis\thanks{\noindent 
Research supported in part by NSF grant DMS-9977354.
{\it Email address:} rellis@math.tamu.edu.  }\\ \\
Department of Mathematics, Texas A\&M University, \\
College Station, TX 77843-3368, USA}
\maketitle

\begin{abstract}
Discrete Green's functions are the inverses or pseudo-inverses of
combinatorial Laplacians.  We present compact formulas for discrete 
Green's functions, in terms of the eigensystems of corresponding Laplacians,
for products of regular graphs with or without boundary.  Explicit
formulas are derived for the cycle, torus, and 3-dimensional
torus, as is an inductive formula for the $t$-dimensional torus
with $n$ vertices, from which the Green's function can be
completely determined in time $O(t\, n^{2-1/t}\log{n})$. These
Green's functions may be used in conjunction with diffusion-like
problems on graphs such as electric potential, random walks, and
chip-firing games or other balancing games.
\end{abstract}


\section{Introduction}
Discrete Green's functions are the inverses or pseudo-inverses of
combinatorial Laplacians, which govern diffusion-like problems on
graphs such as random walks, electric potential and chip-firing
games or other balancing games. Just as Green's functions in the
continuous case depend on the domain and boundary conditions,
discrete Green's functions are associated with the underlying
graph and boundary conditions, if any.  Just as in the continuous
case, a new set of discrete Green's functions must be determined
for each new class of graphs.  Certainly, a discrete Green's
function can be determined by brute force (pseudo-)inversion of
the corresponding Laplacian, but this is no advancement toward
compact or closed-form functions.

In this paper, we develop such compact formulas for discrete
Green's functions on products of simple graphs with or without
boundary. 
Section \ref{sec:prelim} presents the necessary definitions and background.
Section \ref{sec:cycle} illustrates these definitions by deriving the Green's
function for the cycle. In Section \ref{sec:productFormula}, we derive formulas for
products of regular graphs with or without boundary.  In Section
\ref{sec:toriFormula}, we illustrate the case of products of regular graphs without
boundary by deriving the Green's function for the $t$-dimensional
torus and addressing its computational complexity, with explicit
Green's functions given for $t=2$ and $t=3$. Finally, the relationship between
Green's functions and the hitting time of a random walk is explained and 
illustrated in the case of the 2-torus in Section \ref{sec:hittingTime}

\section{Preliminaries\label{sec:prelim}}

The basic definitions follow those of \cite{C97}. Let
$\Gamma=(V,E)$ be a simple connected graph. Let $x,y\in V$ be
arbitrary vertices. The \emph{adjacency matrix} $A$ is defined by
$A(x,y)=\chi(x\sim y)$. The diagonal \emph{degree matrix} $D$ is
defined by $D(x,y)=\chi(x=y)\cdot d_x$, where $d_x$ is the degree
of $x$ in $\Gamma$. The \emph{volume} of a graph is
$vol(\Gamma)=\sum_{v\in V} d_v$. The \emph{Laplacian}, $L=D-A$, of
$\Gamma$ is
$$  L(x,y) = \left\{ \begin{array}{cl}
        d_x, & \mbox{if $x=y$}        \\
        -1, & \mbox{if $x\sim y$}   \\
        0, & \mbox{otherwise.}  \end{array}\right.
$$
The \emph{normalized Laplacian}, $\mathcal{L}=D^{-1/2}LD^{-1/2}$,
is
$$  \mathcal{L}(x,y) = \left\{ \begin{array}{cl}
        1, & \mbox{if $x=y$}        \\
        -1/\sqrt{d_xd_y}, & \mbox{if $x\sim y$} \\
        0, & \mbox{otherwise.}  \end{array}\right.
$$
The \emph{discrete Laplace} operator $\Delta$ is
$$  \Delta(x,y) = \left\{ \begin{array}{cl}
        1, & \mbox{if $x=y$}        \\
        -1/d_x, & \mbox{if $x\sim y$}   \\
        0, & \mbox{otherwise.}  \end{array}\right.
$$
In general, the following relations hold between $L$,
$\mathcal{L}$, and $\Delta$:
\begin{eqnarray}
    \mathcal{L} & = & D^{-1/2}LD^{-1/2} \nonumber \\
    \mathcal{L} & = & D^{1/2}\Delta D^{-1/2} \nonumber \\
    L & = & D\Delta.    \nonumber
\end{eqnarray}

For $S\subseteq V$, we define the \emph{Dirichlet} versions of
$L_S$, $\mathcal{L}_S$, and $\Delta_S$ as the results of deleting
the rows and columns corresponding to $V\setminus S$ from $L$,
$\mathcal{L}$, and $\Delta$, respectively.  Without loss of generality,
we assume that both $\Gamma$ and the subgraph generated by $S$ to be
connected. By abuse of notation, we use $S$ to refer to the subgraph
generated by $S$ in $\Gamma$.  Additionally, we say that
the orthonormal eigensystem of $S$ is
the orthonormal eigensystem $\{(\lambda_j,\phi_j):j\in J\}$ of the
real symmetric matrix $\mathcal{L}_S$;
where $J=\{0,\ldots,|V|-1\}$ with $\lambda_0=0$ when $S=V$, and 
$J=\{1,\ldots,|S|\}$ with $\lambda_1>0$ when $S\neq V$.  In either case, 
we order the eigenvalues by $\lambda_j\leq \lambda_{j+1}$ for all 
subscripts in range; basic properties such as $0\leq \lambda_j\leq 2$ for all $j$ 
are summarized in \cite{CY00}.  In particular, the orthonormal eigenvectors
can be chosen to have real entries.

When $S\subsetneq V$, $\Delta_S$, $L_S$ and $\mathcal{L}_S$ are
invertible, and the \emph{Green's function} $G$ and
\emph{normalized Green's function} $\mathcal{G}$ are determined by
the relations
\begin{eqnarray}
& \Delta_S G = G \Delta_S =  I_S & \qquad\mbox{and}  \nonumber \\
& \mathcal{L}_S\mathcal{G} = \mathcal{G}\mathcal{L}_S  = I_S. &
    \label{eqn:normalizedGMatrixRelation}
\end{eqnarray}
We can tie these relations to random walks as follows.  Let
$P=[P(x,y)]$ be the transition probability matrix for the simple
irreducible \emph{transient} random walk on $S$ with absorbing
states $V\setminus S$, where the probability $p_{xy}$ of moving to
state $y$ from state $x$ is $1/d_x$ if $x$ and $y$ are adjacent
and 0 otherwise. Then $\Delta_S=I-P$, and
$(I-P)^{-1}=I+P+P^2+\cdots$ gives
\begin{equation}\label{eqn:transientRandomWalkG}
G(x,y) = \sum_n P_n(x,y),
\end{equation}
where $P_n(x,y)$ is the $n$-step transition probability matrix
(cf. \cite[p.\ 31]{T00}). See \cite{AF02, DS84} for definitions
and results on random walks.

On the other hand, when $S=V$, the sum
in (\ref{eqn:transientRandomWalkG}) does not converge, and we
require an alternate definition of the Green's function.  Since
$\mathcal{L}$ is not invertible, its corresponding normalized
Green's function $\mathcal{G}$ is defined by
\begin{equation}
\mathcal {G} = \sum_{\lambda_j>0}\frac{1}{\lambda_j}\phi_j
\phi_j^*. \label{eqn:normalizedGFourier}
\end{equation}
This definition of $\mathcal{G}$ is equivalent to the two relations
\begin{eqnarray}
& \mathcal{G}\mathcal{L} = \mathcal{L}\mathcal{G}  =
    I - P_0 = I-\phi_0\phi_0^* & \qquad\mbox{and} \nonumber\\
& \mathcal{G}P_0 = 0, &
\label{eqn:normalizedGPseudoInverseRelations}
\end{eqnarray}
where $P_0=\phi_0\phi_0^*$ is the projection of the orthonormal
eigenvector $\phi_0$ corresponding to eigenvalue 0.  Since
$L\vec{1}=0$, we have $\phi_0=D^{1/2}\vec{1}/||D^{1/2}\vec{1}||$,
and so $\phi_0(x)=\sqrt{d_x/vol(\Gamma)}$. It is also important to
note that when $\mathcal{G}$ is invertible, the definitions in
(\ref{eqn:normalizedGMatrixRelation}) and
(\ref{eqn:normalizedGFourier}) are equivalent.  Furthermore,
$\mathcal{G}$ is related to the so-called {\em fundamental
matrix\/}
\begin{equation}\label{eqn:fundMatFormal} 
Z(x,y):=\sum_{n=0}^\infty(P_n(x,y)-\pi_y),
\end{equation}
where $\pi$ is the stationary distribution of the random walk on
$\Gamma$, by the equation $\mathcal{G}=D^{1/2}ZD^{-1/2}$; simply
verify that $D^{1/2}ZD^{-1/2}$ satisfies
(\ref{eqn:normalizedGFourier}) and
(\ref{eqn:normalizedGPseudoInverseRelations}). See Chapter 3, p.\
17 of \cite{AF02} for relationships between $Z$ and hitting times.
A more complete exposition on discrete Green's functions in the
context of spectral graph theory appears in \cite{E02}.

\section{Green's function for the cycle\label{sec:cycle}}
In this section we illustrate the preceding definitions
in the case of the cycle, which has no boundary.  The techniques
developed here will be used in later sections to construct the
Green's functions for higher dimensional tori.

By the cycle $C_m$, we mean the 2-regular connected graph with vertices
$\{0,1,\ldots,m-1\}$, where $m\geq 3$.
The various Laplacians are related by
$\Delta=\mathcal{L}=L/2$. Applying the definition in
(\ref{eqn:normalizedGPseudoInverseRelations}), the normalized
Green's function $\mathcal{G}$ satisfies
\begin{eqnarray}
\mathcal{G}\mathcal{L} \ = \ \mathcal{L}\mathcal{G} & = &
    I - \frac{1}{m}J \qquad\mbox{and} \nonumber\\
\mathcal{G}J &  =  & 0,
    \label{eqn:normalizedGCyclePseudoInverseRelations}
\end{eqnarray}
where $\phi_0(x)=\sqrt{1/m}$, and $J$ is the $m\times m$ matrix
of 1's.  Because the cycle is vertex-transitive, 
the values $\mathcal{L}(x,y)$ and
$\mathcal{G}(x,y)$ depend only on the distance $|y-x|$ between $x$
and $y$, and the following definition of $\mathcal{G}(a)$ is
well-defined:
\begin{equation}
\mathcal{G}(a) := \mathcal{G}(x,y), \qquad\mbox{if $a=|y-x|$.}
 \label{eqn:normalizedGCyclePseudoInverseRelationsA}
\end{equation}
We are ready to derive the Green's function for the cycle.

\begin{theorem}\label{thm:normalizedGCycle}
Let $m\geq 3$.  For \ $0\leq x,y\leq m-1$, the cycle $C_m$ has
normalized Green's function
\begin{equation}\label{eqn:normalizedGCycle}
\mathcal{G}(x,y) = \frac{(m+1)(m-1)}{6m} - |y-x| +
    \frac{(y-x)^2}{m}.
\end{equation}
\end{theorem}
\begin{proof} From (\ref{eqn:normalizedGCyclePseudoInverseRelations}) and
(\ref{eqn:normalizedGCyclePseudoInverseRelationsA}), we
have the recurrence
\begin{eqnarray}
2\mathcal{G}(x,y)-\mathcal{G}(x,y-1)-\mathcal{G}(x,y+1) & = &
    \left\{\begin{array}{cl}
        2-2/m, & x=y    \\
        -2/m, & x\neq y,\end{array}\right.\quad \mbox{or}
    \nonumber \\
2\mathcal{G}(a)-\mathcal{G}(a-1)-\mathcal{G}(a+1) & = &
    \left\{\begin{array}{cl}
        2-2/m, & a=0,    \\
        -2/m, & a > 0,\end{array}\right.
    \nonumber
\end{eqnarray}
provided that we define $\mathcal{G}(-1)=\mathcal{G}(1)$ for
simplicity of representing the case $a=0$. The following
recurrence on differences results:
\begin{eqnarray}
    \mathcal{G}(a+1)-\mathcal{G}(a) & = & \mathcal{G}(a)
    -\mathcal{G}(a-1) \  +  \ \frac{2}{m} \ - \ 2\chi(a=0).
    \label{eqn:normalizedGCycleDiffRecur}
\end{eqnarray}
The second constraint in
(\ref{eqn:normalizedGCyclePseudoInverseRelations}) determines that
the sum of $\mathcal{G}$ across any row must be 0; i.e.,
\begin{eqnarray}
    \sum_{a=0}^{m-1}\mathcal{G}(a) & = & 0. \label{eqn:normalizedGCycleRowsum}
\end{eqnarray}
We now solve the recurrence on differences, starting with
$\mathcal{G}(1)-\mathcal{G}(0)$ by setting $a=0$ in
(\ref{eqn:normalizedGCycleDiffRecur}).
\begin{eqnarray}
    \frac{1}{2}\left(\mathcal{G}(1)
    -2\mathcal{G}(0)+\mathcal{G}(-1)\right)
    & = & \mathcal{G}(1)-\mathcal{G}(0) \nonumber \\
    & = & \frac{1}{m} - 1.      \label{eqn:normalizedGCycleDiffBaseCase}
\end{eqnarray}
Resolving (\ref{eqn:normalizedGCycleDiffRecur}) with base case
given by (\ref{eqn:normalizedGCycleDiffBaseCase}), we have
\begin{eqnarray}
\mathcal{G}(a+1)-\mathcal{G}(a) & = &
    \mathcal{G}(a)-\mathcal{G}(a-1)+\frac{2}{m} \nonumber \\
& = & \mathcal{G}(a-1)-\mathcal{G}(a-2)+2 \cdot \frac{2}{m} \nonumber \\
& = & \quad \vdots \nonumber \\
& = & \mathcal{G}(1)-\mathcal{G}(0)+a \cdot \frac{2}{m} \nonumber \\
& = & \frac{1}{m}-1+a \cdot \frac{2}{m} \ .
\label{eqn:normlizedGCycleRecur}
\end{eqnarray}

Having derived a simple recurrence from the recurrence on differences, 
we proceed to determine $\mathcal{G}(a)$. Resolving
(\ref{eqn:normlizedGCycleRecur}) yields
\begin{eqnarray}
\mathcal{G}(a) & = & \mathcal{G}(a-1) + \frac{1}{m}-1+(a-1)
    \frac{2}{m} \nonumber \\
& = & \mathcal{G}(a-2) + \frac{2}{m}-2 + [(a-1)+(a-2)]\frac{2}{m}
    \nonumber \\
& = & \quad \vdots \nonumber \\
& = & \mathcal{G}(0) + \frac{a}{m}-a+\frac{2}{m}
    \sum_{k=0}^{a-1}k   \nonumber \\
& = & \mathcal{G}(0) - a + \frac{a^2}{m} \ .
    \label{eqn:normalizedGCycleAG0}
\end{eqnarray}
Now applying the row sum constraint
(\ref{eqn:normalizedGCycleRowsum}) allows us to compute the value
of $\mathcal{G}(0)$.
\begin{eqnarray}
0 & = & \sum_{a=0}^{m-1}
    \left(\mathcal{G}(0) - a + \frac{a^2}{m}\right) \nonumber \\
m \cdot \mathcal{G}(0) & = & \sum_{a=0}^{m-1}
    \left(a-\frac{a^2}{m}\right) \nonumber\\
\mathcal{G}(0) & = & \frac{(m+1)(m-1)}{6m} \ .
    \label{eqn:normalizedGCycleG0}
\end{eqnarray}
Plugging (\ref{eqn:normalizedGCycleG0}) into
(\ref{eqn:normalizedGCycleAG0}) and letting $a=|y-x|$ achieves the
desired result. 
\end{proof}

Because (\ref{eqn:normalizedGFourier}) also gives $\mathcal{G}$
for $C_m$, we have a whole class of identities formed by choosing
any orthonormal eigenbasis for $C_m$ and equating
(\ref{eqn:normalizedGFourier}) with (\ref{eqn:normalizedGCycle}).
We give one such well-known basis now, which arises naturally from the
consideration of circulant matrices (cf.\ \cite{D79}).

\begin{lemma}\label{lem:cycleMixedEigensystem}
For $m\geq 3$, define 
$\phi_j(x):=\frac{1}{\sqrt{m}}\exp{\left(-i\frac{2\pi
jx}{m}\right)}$.
Then $\{\big(1-\cos{(2\pi
j/m)},\phi_j\big):0\leq j<m\}$ is an orthonormal eigensystem of
$\mathcal{L}$ for $C_m$.
\end{lemma}
The proof is omitted, but follows from a straightforward
verification of the necessary conditions.
Theorem \ref{thm:cycleGreenIdentity}
follows by combining Theorem
\ref{thm:normalizedGCycle} with (\ref{eqn:normalizedGFourier})
using the orthonormal eigenbasis of Lemma
\ref{lem:cycleMixedEigensystem}. 

\begin{theorem}\label{thm:cycleGreenIdentity} Let $m\geq 3$ and let $0\leq x,y<m$.  Then
$$\frac{1}{m}\sum_{j=1}^{m-1}
    \frac{\exp{\big((2\pi ij/m)(y-x)\big)}}{1-\cos{(2\pi j/m)}}
    = \frac{(m+1)(m-1)}{6m} - |y-x| +
    \frac{(y-x)^2}{m} \ .
$$
\end{theorem}

\section{Green's functions for products of regular graphs
\label{sec:productFormula}}

Theorems 4-5 of \cite{CY00} give a contour integral formula 
for the Green's function of the Cartesian product of two regular graphs
with boundary, provided a certain generalized Green's function is 
known for each factor graph. In this section, we extend these 
results to include the cases where one or both of the factor graphs 
is without boundary, and provide simplified working formulas requiring 
the generalized Green's function of one graph and the eigensystem of
the other. 
The original technique of \cite{CY00} recovers
the Green's function of the product graph as the residues of a certain 
contour 
integral whose contour contains as poles all of the eigenvalues of one of the graphs,
and none of the negatives of the eigenvalues of the other graph.  
When both 
graphs are without boundary, both have an eigenvalue of 0, and 
the resulting  order 2 pole requires a modified contour integral formula.

Recalling the definitions in Section \ref{sec:prelim},
Let $\Gamma=(V,E)$ and $\Gamma'=(V',E')$ be regular graphs of degree 
$d$ and $d'$ with specified
vertex subsets $S$ and $S'$ and normalized Dirichlet Laplacians
$\mathcal{L}_S$ and $\mathcal{L}'_{S'}$, respectively.  
For any $\alpha\in\mathbb{C}$, define $\mathcal{G}_\alpha$ to be the
symmetric matrix satisfying the relation
$(\mathcal{L}_S+\alpha)\mathcal{G}_\alpha=I_S$, if $S\subsetneq
V$; and the relations
\begin{eqnarray}
(\mathcal{L}_S+\alpha)\mathcal{G}_\alpha & = & I_S-P_0,
    \qquad\mbox{and} \nonumber \\
\mathcal{G}_\alpha P_0&=&0,
    \label{eqn:normalizedGAlphaPseudoInverseRelations}
\end{eqnarray}
if $S=V$ ($\mathcal{L}_S$ is singular). 
In either case, this is equivalent to
\begin{equation}
\mathcal{G}_{\alpha}(x,y) = \sum_{\lambda_j >
0}\frac{1}{\lambda_j+\alpha}\phi_j(x)\overline{\phi_j(y)},
\label{eqn:normalizedGAlphaFourier}
\end{equation}
where the $\phi_j$'s are the orthonormal eigenfunctions of
$\mathcal{L}_S$ associated with the eigenvalues $\lambda_j$. In particular,
$\mathcal{G}_\alpha$ is a rational function of $\alpha$. The
analogous definitions  of
$\mathcal{G}'_\alpha$, $\phi'_i$, and $\lambda'_i$ are made for
$\Gamma'$.

The choice of $S$ and $S'$ induces the specified vertex set $S\times S'$
in the Cartesian product $\Gamma\times\Gamma'$.  This product 
of $\Gamma=(V,E)$ and $\Gamma'=(V',E')$ has vertex set $\{(v,v'):v\in V,v'\in
V'\}$ and edges of the form $\{(v,v'),(v,u')\}$ or
$\{(v,v'),(u,v')\}$ where $\{u,v\}\in E$, $\{u',v'\}\in E'$. Let $\mathcal{L}^{\times}$
be the normalized Laplacian of  $\Gamma\times\Gamma'$, with corresponding 
normalized Dirichlet Laplacian $\mathcal{L}^{\times}_{S\times S'}$.  
We abuse notation
by referring to $S$, $S'$ and $S\times S'$ as the graphs they induce.
If $\mathcal{L}_S$ and $\mathcal{L}'_{S'}$ have orthonormal eigensystems 
$\{(\lambda_j,\phi_j):j\in J\}$ and $\{(\lambda'_k,\phi'_k):k\in K\}$, 
respectively, then when $S$ and $S'$ are regular of the same degree, 
$\mathcal{L}^{\times}_{S\times S'}$ has orthonormal eigensystem 
\begin{equation}\label{eqn:productEigensystem}
\left\{\left(\frac{\lambda_j+\lambda'_k}{2},\Phi_{j,k}\right)
:j\in J, k\in K\right\},
\end{equation}
where $\Phi_{j,k}(v,v'):=\phi_j(v)\phi'_k(v').$
We begin with the case in which $S$ has boundary.  

\subsection{At least one graph has boundary\label{sec:resultsBoundary}}

Suppose $S\subsetneq V$, so that $S$ generates a connected
subgraph with boundary in $\Gamma$.  We allow $S'\subseteq V'$ to
generate an arbitrary connected subgraph in $\Gamma'$. 
We consider two cases, where the factor graphs are regular of the 
same degree or regular of different degrees.

First, suppose $S$ and $S'$ be regular of the same degree, and 
let $C$ denote a simple closed contour in the complex plane 
consisting of all
$\alpha\in\mathbb{C}$ satisfying $|\alpha-1|=1+\lambda_1/2$. The
contour $C$ is designed to enclose all of the $\lambda'_k$'s and
none of the $-\lambda_j$'s.  We have the following generalization
of Theorem 4 of \cite{CY00}.

\begin{theorem}\label{thm:normalizedGProductBoundary}
Let $S$ and $S'$ be induced subgraphs of\, $\Gamma=(V,E)$ and
$\Gamma'=(V',E')$, respectively, which are both regular of degree
$d$. Let $S\subsetneq V$ and $S'\subseteq V'$.  The normalized
Green's function $\mathbf{G}$ of the Cartesian product $S\times
S'$ with Dirichlet boundary condition is
$$ \mathbf{G}((x,x'),(y,y'))=\frac{1}{\pi
i}\int_C\mathcal{G}_\alpha(x,y)\mathcal{G}'_{-\alpha}(x',y')d\alpha.
$$
\end{theorem}
\begin{proof} Combining the product graph eigensystem in
(\ref{eqn:productEigensystem}) with the 
formal definition of $\mathbf{G}$ in
(\ref{eqn:normalizedGFourier}), we have
\begin{eqnarray}
\mathbf{G}((x,x'),(y,y')) & = &
    2\sum_{j,k}\frac{\Phi_{j,k}(x,x')\overline{\Phi_{j,k}(y,y')}}
    {\lambda_j+\lambda'_k} \nonumber \\
& = & \frac{1}{\pi i}\int_C \sum_{j=1}^{|S|}\sum_k
    \frac{\phi_j(x)\overline{\phi_j(y)}\phi'_k(x')
    \overline{\phi'_k(y')}}{(\lambda_j+\alpha)
    (\lambda'_k-\alpha)}d\alpha \nonumber \\
& = & \frac{1}{\pi i}\int_C
    \mathcal{G}_\alpha(x,y)\mathcal{G}'_{-\alpha}(x',y')d\alpha,
    \qquad\mbox{by (\ref{eqn:normalizedGAlphaFourier}).} \nonumber
    \qquad\qquad
\end{eqnarray}
\end{proof}

Note that the above contour integral picks up exactly the residues
at $\alpha=\lambda'_k$.  For example, the residue at
$\alpha=\lambda'_{k}$ is exactly
\begin{equation}
\sum_{\lambda'_K,\lambda'_{K}=\lambda'_k}
\sum_{j=1}^{|S|}\frac{\phi_j(x)\overline{\phi_j(y})
    \phi'_K(x')\overline{\phi'_K(y')}}
    {(\lambda_j+\lambda'_{k})}=
\sum_{\lambda'_K,\lambda'_{K}=\lambda'_k}\phi'_{K}(x')
    \overline{\phi'_{K}(y')}
    \cdot\mathcal{G}_{\lambda'_{k}}(x,y).
\label{eqn:residueGProductBoundary}
\end{equation}
For convenience, we assign the term of the residue in
(\ref{eqn:residueGProductBoundary}) corresponding to $K$ to
$\lambda'_K$. This observation gives us the computational formula
for $\mathbf{G}$ in the following corollary, which in practice may be 
applied to yield a closed formula for a specific product graph, or
to generate in conjunction with (\ref{eqn:normalizedGFourier}) 
a non-trivial identity involving 
$\mathcal{G}_\alpha$ and arbitrary orthonormal eigensystems of $\mathcal{L}_S$, 
$\mathcal{L}_{S'}$ and $\mathcal{L}^{\times}_{S\times S'}$.

\begin{cor}\label{cor:normalizedGProductBoundary}
Under the same conditions as in Theorem
\ref{thm:normalizedGProductBoundary}, we have
$$\mathbf{G}((x,x'),(y,y')) = 2\sum_k
\phi'_k(x')\overline{\phi'_k(y')}\mathcal{G}_{\lambda'_k}(x,y).$$
\end{cor}

Second, suppose we have the same conditions as Theorem
\ref{thm:normalizedGProductBoundary}, except $S$ and $S'$
are regular of degrees $d$ and $d'$,
respectively.  The normalized Dirichlet Laplacian 
$\mathcal{L}^{\times}_{S\times S'}$ has orthonormal 
eigensystem
\begin{equation}\label{eqn:productEigensystem2}
\left\{\left(\frac{d}{d+d'}\lambda_j+\frac{d'}{d+d'}\lambda'_k,\Phi_{j,k}\right)
:j\in J, k\in K\right\}.
\end{equation}
The poles
of $\mathcal{G}_{\alpha/d}$ are at $\alpha=-d\lambda_j$, and the
poles of $\mathcal{G}_{-\alpha/d'}$ are at $\alpha=d'\lambda'_k$.
Let $C$ denote a simple closed contour in the complex plane consisting of all
$\alpha\in\mathbb{C}$ satisfying $|\alpha-d'|=d'+d\lambda_1/2$;
thus $C$ contains all of the $d'\lambda'_k$'s but none of the
$-d\lambda_j$'s.  We obtain the following minor extension to
Theorem 5 of \cite{CY00}.
\begin{theorem}
\label{thm:normalizedGProductBoundaryDDPrime} Let $S$ and $S'$ be
induced subgraphs of $\Gamma=(V,E)$ and $\Gamma'=(V',E')$,
respectively, where $\Gamma$ is regular of degree $d$ and
$\Gamma'$ is regular of degree $d'$. Let $S\subsetneq V$ and
$S'\subseteq V'$. The normalized Green's function $\mathbf{G}$ of
the Cartesian product $S\times S'$ with Dirichlet boundary
condition is
$$ \mathbf{G}((x,x'),(y,y'))=\frac{d+d'}{2\pi idd'}
\int_C\mathcal{G}_{\alpha/d}(x,y)\mathcal{G}'_{-\alpha/d'}(x',y')d\alpha.
$$
\end{theorem}
\begin{proof} Combining (\ref{eqn:productEigensystem2})
with (\ref{eqn:normalizedGFourier}), we have
\begin{eqnarray}
\mathbf{G}((x,x'),(y,y')) & = & \sum_{j,k}
    \frac{d+d'}{d\lambda_j+d'\lambda'_k}
    \Phi_{j,k}(x,x')\overline{\Phi_{j,k}(y,y')} \nonumber\\
& = & \frac{d+d'}{2\pi i}\int_C
    \sum_{j=1}^{|S|}\sum_k\frac{\phi_j(x)\overline{\phi_j(y)}
    \phi'_k(x')\overline{\phi'_k(y')}}
    {(d\lambda_j+\alpha)(d'\lambda'_k-\alpha)}    \nonumber \\
& = & \frac{d+d'}{2\pi idd'}\int_C
    \sum_{j=1}^{|S|}\sum_k\frac{\phi_j(x)\overline{\phi_j(y)}
    \phi'_k(x')\overline{\phi'_k(y')}}
    {(\lambda_j+\alpha/d)(\lambda'_k-\alpha/d')}    \nonumber \\
& = & \frac{d+d'}{2\pi idd'}\int_C \mathcal{G}_{\alpha/d}(x,y)
    \mathcal{G}'_{-\alpha/d'}(x',y')d\alpha. \nonumber
    \qquad\qquad
\end{eqnarray}
\end{proof}

Analogous to Corollary \ref{cor:normalizedGProductBoundary}, by
inspecting the residues of the above contour integral at all
values $\alpha=d'\lambda'_k$, we have the following.
\begin{cor} Under the same conditions as Theorem
\ref{thm:normalizedGProductBoundaryDDPrime}, we have
$$
\mathbf{G}((x,x'),(y,y')) = \frac{d+d'}{d}\sum_k
\phi'_k(x')\overline{\phi'_k(y')}
\mathcal{G}_{d'\lambda'_k/d}(x,y).\nonumber
$$
\end{cor}

\subsection{Neither graph has boundary\label{sec:resultsNoBoundary}}

Here we consider the case of $S=V$ and $S'=V'$, so that $S\times S'$
is the entire product graph $\Gamma\times\Gamma'$; recall that the
normalized Green's function for $\Gamma\times\Gamma'$ is 
non-invertible. Let $m=|V|$ and
$n=|V'|$, and consider first the case in which $\Gamma$ and $\Gamma'$ 
have the same degree. 

The orthonormal eigensystems of $\mathcal{L}$ and $\mathcal{L}'$ are
$\{(\phi_j,\lambda_j):0\leq j\leq m-1\}$ and
$\{(\phi'_k,\lambda'_k):0\leq k\leq n-1\}$, respectively. 
Let $C$ denote a simple closed contour
in the complex plane consisting of all $\alpha\in\mathbb{C}$
satisfying $|\alpha-(2+\lambda'_1/2)|=2$. The contour $C$ is
designed to enclose $\lambda'_1,\ldots,\lambda'_{n-1}$ and none of
$-\lambda_0,\ldots,-\lambda_{m-1}$ or $\lambda'_0$.

\begin{theorem} \label{thm:normalizedGProductNoBoundary}
Let\/ $\Gamma$ and $\Gamma'$ be connected regular graphs of
degrees $d$ without boundary. With the notation above, the
normalized Green's function $\mathbf{G}$ of the Cartesian product
$\Gamma\times\Gamma'$ is
\begin{eqnarray}
\mathbf{G}((x,x'),(y,y')) & = & \frac{1}{\pi i} \int_C
    \mathcal{G}_\alpha(x,y)\mathcal{G}'_{-\alpha}(x',y')d\alpha
    + \frac{2}{n}
    \mathcal{G}(x,y)
    +\frac{2}{m}
    \mathcal{G}'(x',y'). \nonumber
\end{eqnarray}
\end{theorem}
\begin{proof} The eigenvector $\phi_0$ corresponding to eigenvalue 0 is
determined by $\phi_0(x)=\sqrt{(d_x/vol(\Gamma))}$. Combining 
the formal definition of $\mathbf{G}$ in (\ref{eqn:normalizedGFourier})
with (\ref{eqn:productEigensystem}), and noting that $d_v=d_{v'}=d$,
$vol(\Lambda)=d\cdot m$, and $vol(\Lambda')=d\cdot n$, we have
\begin{eqnarray}
    \mathbf{G}((x,x'),(y,y')) & = &
    2\sum_{j=1}^{m-1}\sum_{k=1}^{n-1}
    \frac{\Phi_{j,k}(x,x')\overline{\Phi_{j,k}(y,y')}}
    {\lambda_j+\lambda'_k}\nonumber \\
& & + 2\sum_{j=1}^{m-1}\frac{\Phi_{j,0}(x,x')
    \overline{\Phi_{j,0}(y,y')}} {\lambda_j}
    + 2\sum_{k=1}^{n-1}
    \frac{\Phi_{0,k}(x,x')\overline{\Phi_{0,k}(y,y')}}
    {\lambda'_k} \nonumber\\
& = & \frac{1}{\pi i}\int_C \sum_{j=1}^{m-1}\sum_{k=1}^{n-1}
    \frac{\phi_j(x)\overline{\phi_j(y)}
    \phi_k(x')\overline{\phi_k(y')}}{(\lambda_j+\alpha)
    (\lambda'_k-\alpha)}d\alpha \nonumber \\
& & + 2\frac{\sqrt{d_{x'}d_{y'}}}{vol(\Gamma')}
    \sum_{j=1}^{m-1}\frac{\phi_j(x)\overline{\phi_j(y)}
    }{\lambda_j }
    + 2\frac{\sqrt{d_{x}d_{y}}}{vol(\Gamma)}
    \sum_{k=1}^{n-1}\frac{\phi_k(x')\overline{\phi_k(y')}}{\lambda'_k}
    \nonumber\\
& = & \frac{1}{\pi i}\int_C
    \mathcal{G}_\alpha(x,y)\mathcal{G}'_{-\alpha}(x',y')d\alpha
+ \frac{2}{n}
    \mathcal{G}(x,y)
    +\frac{2}{m}
    \mathcal{G}'(x',y')
    \qquad\mbox{by (\ref{eqn:normalizedGAlphaFourier}).} \nonumber
    \qquad
\end{eqnarray}
\end{proof}

Analogous to Corollary \ref{cor:normalizedGProductBoundary},
inspecting the residues of the above contour integral at
$\lambda'_1,\ldots,\lambda'_{n-1}$ yields the following corollary.

\begin{cor}\label{cor:normalizedGProductWOBoundary}
Under the same conditions as in Theorem
\ref{thm:normalizedGProductNoBoundary}, we have
$$\mathbf{G}((x,x'),(y,y')) = 2\sum_{k=1}^{n-1}
\phi'_k(x')\overline{\phi'_k(y')} \mathcal{G}_{\lambda'_k}(x,y) +
\frac{2}{n} \mathcal{G}(x,y) +\frac{2}{m} \mathcal{G}'(x',y'). $$
\end{cor}

Now suppose we have the same conditions as Theorem
\ref{cor:normalizedGProductWOBoundary}, except that the graphs
$\Gamma$ and $\Gamma'$ are regular of degrees $d$ and $d'$,
respectively, and  $\mathcal{L}^{\times}_{\Gamma\times\Gamma'}$ has
orthonormal eigensystem given by (\ref{eqn:productEigensystem2}).
Let $C$ denote a contour in the
complex plane consisting of all $\alpha\in\mathbb{C}$ satisfying
$|\alpha-(d'+d'\lambda'_1/2)|=d'$. Thus $C$ contains
$d'\lambda'_1,\cdots,d'\lambda'_{n-1}$, but neither
$-d\lambda_0,\ldots, -d\lambda_{m-1}$ nor $d'\lambda'_0$. We obtain
the following theorem.

\begin{theorem}
\label{thm:normalizedGProductNoBoundaryDDPrime} Let\/ $\Gamma$ and
$\Gamma'$ be connected regular graphs without boundary of degrees
$d$ and $d'$, respectively. With the notation above, the normalized Green's
function $\mathbf{G}$ of the Cartesian product
$\Gamma\times\Gamma'$ is
$$ \mathbf{G}((x,x'),(y,y'))=\frac{d+d'}{2\pi idd'}
\int_C\mathcal{G}_{\alpha/d}(x,y)\mathcal{G}'_{-\alpha/d'}(x',y')d\alpha
+ \frac{d+d'}{dn}\mathcal{G}(x,y)
    + \frac{d+d'}{d'm}\mathcal{G}'(x',y').
$$
\end{theorem}
\begin{proof} Combining the definition of $\mathbf{G}$ in
(\ref{eqn:normalizedGFourier}) with (\ref{eqn:productEigensystem2}), we have
\begin{eqnarray}
\mathbf{G}((x,x'),(y,y'))
& = & \sum_{j=1}^{m-1}\sum_{k=1}^{n-1}
    \frac{d+d'}{d\lambda_j+d'\lambda'_k}
    \Phi_{j,k}(x,x')\overline{\Phi_{j,k}(y,y')} \nonumber\\
& & +
    \sum_{j=1}^{m-1}\frac{d+d'}{d\lambda_j}\Phi_{j,0}(x,x')
    \overline{\Phi_{j,0}(y,y')}
    + \sum_{k=1}^{n-1}\frac{d+d'}{d'\lambda'_k}\Phi_{0,k}(x,x')
    \overline{\Phi_{0,k}(y,y')} \nonumber\\
& = & \frac{d+d'}{2\pi idd'}
    \int_C\frac{dd'\phi_j(x)\overline{\phi_j(y)}
    \phi_k(x')\overline{\phi_k(y')}}
    {(d\lambda_j+\alpha)(d'\lambda_k-\alpha)}d\alpha
    \nonumber\\
& & + \frac{d+d'}{d}\frac{\sqrt{d_{x'}d_{y'}}}{vol(\Gamma')}
    \mathcal{G}(x,y)
    + \frac{d+d'}{d}\frac{\sqrt{d_xd_y}}{vol(\Gamma)}
    \mathcal{G'}(x',y') \nonumber \\
& = & \frac{d+d'}{2\pi idd'}\int_C \mathcal{G}_{\alpha/d}(x,y)
    \mathcal{G}'_{-\alpha/d'}(x',y')d\alpha
    + \frac{d+d'}{dn}\mathcal{G}(x,y)
    + \frac{d+d'}{d'm}\mathcal{G}'(x',y').
     \nonumber
\end{eqnarray}
\end{proof}

Analogous to Corollary \ref{cor:normalizedGProductBoundary}, by
inspecting the residues of the contour integral at
$d'\lambda'_1,\ldots,d'\lambda'_{n-1}$, we have the following.

\begin{cor}\label{cor:normalizedGProductWOBoundaryDDPrime}
Under the same conditions as in Theorem
\ref{thm:normalizedGProductNoBoundaryDDPrime}, we have
$$\mathbf{G}((x,x'),(y,y')) = \frac{d+d'}{d}\sum_{k=1}^{n-1}
\phi'_k(x')\overline{\phi'_k(y')}
\mathcal{G}_{d'\lambda'_k/d}(x,y) + \frac{d+d'}{dn}
\mathcal{G}(x,y) + \frac{d+d'}{d'm} \mathcal{G}'(x',y'). $$
\end{cor}

\section{Example Green's functions for products of 
cycles\label{sec:toriFormula}}

Application of the results in the previous two sections to
specific examples is limited only to cases in which the necessary raw
materials can be computed.  The results in Section
\ref{sec:resultsBoundary} for the product where at least one graph
has boundary require the eigensystem of one graph and the
generalized Green's function $\mathcal{G}_\alpha$ of the other
graph.  Although as written, the results assume that the
generalized Green's function is known for the graph without
boundary, and the eigensystem is known for the graph without
boundary, Theorems \ref{thm:normalizedGProductBoundary} and
\ref{thm:normalizedGProductBoundaryDDPrime} can be easily
re-derived for when the reverse is true.

The results in Section \ref{sec:resultsNoBoundary} for the product
of two graphs without boundary require knowledge of the
eigensystem of one graph, the generalized Green's function
$\mathcal{G}_\alpha$ of the other graph, and the Green's
functions $\mathcal{G}$ and $\mathcal{G'}$ of both graphs.  
Thus computation of the normalized Green's function for any graph can 
be done via any decomposition of the graph into factor graphs where 
this  information is known. This
observation can be particularly useful in computing the normalized
Green's function inductively where each additional factor graph is
from a specific family.

\subsection{The torus $C_m\times C_n$}

Following Corollary \ref{cor:normalizedGProductWOBoundary},
determination of the Green's function of the torus requires the
Greens function $\mathcal{G}$ and generalized Green's function
$\mathcal{G}_\alpha$ of the cycle. In obtaining a compact formula
for the torus, it is critical to simplify $\mathcal{G}_\alpha$ as
much as possible before incorporating it into Corollary
\ref{cor:normalizedGProductWOBoundary}.

\begin{theorem}\label{thm:normalizedGTorus}
Let $m,n\geq 3$.  For\, $0\leq x,y\leq m-1$ and $0\leq x',y'\leq
n-1$, the torus \ $C_m\times C_n$ has normalized Green's function
\begin{align*}
\mathbf{G}((x,&x'),(y,y'))  =
    \frac{2}{n}\sum_{k=1}^{n-1}
    \exp\big((2\pi ik/n)(y'-x')\big)
    \left[-\frac{1}{m(1-\cos{(2\pi k/n)})}\right. \nonumber \\
& \left.+\frac{T_{m/2-|y-x|}(2-\cos{(2\pi k/n)})}
    {(1-\cos{(2\pi k/n)})
    (3-\cos{(2\pi k/n)})U_{m/2-1}(2-\cos{(2\pi k/n)})}
    \right]\nonumber \\
& +\frac{2}{n}\left(\frac{(m+1)(m-1)}{6m} - |y-x| +
    \frac{(y-x)^2}{m}\right) \nonumber \\
& +\frac{2}{m}\left(\frac{(n+1)(n-1)}{6n} - |y'-x'| +
    \frac{(y'-x')^2}{n}\right) \nonumber
\end{align*}
where  $T$ and $U$ are the Chebyshev polynomials of the first and
second kind, respectively.
\end{theorem}

Note that the formula depends only on the distances between $y$
and $x$ and between $y'$ and $x'$, which is expected due to the
translational symmetries of the torus.  The first step of the
proof is to obtain a closed form for $\mathcal{G}_\alpha$ for the
cycle $C_m$ ($\mathcal{G}$ was determined in Theorem
\ref{thm:normalizedGCycle}). The proof of Theorem
\ref{thm:normalizedGTorus} is deferred until after Cor.\@
\ref{cor:normalizedGAlphaCycle}.

\begin{theorem} \label{thm:normalizedGAlphaCycleR}
Let $m\geq 3$.  For a cycle $C_m$ with vertices $0,1,\ldots,m-1$, complex
$\alpha\neq 0$, and $0\leq x, y\leq m-1$,  the generalized Green's
function $\mathcal{G}_\alpha$ satisfies
\begin{equation*}
\mathcal{G}_\alpha(x,y) \ = \ -\frac{2}{m\left(r+r^{-1}-2\right)}
    + \frac{2(r^{m/2-|x-y|} +
    r^{-m/2+|x-y|})} {(r-r^{-1})(r^{m/2}-r^{-m/2})},
\end{equation*}
where $2(1+\alpha)=r+r^{-1}$.
\end{theorem}
\begin{proof}
Because $C_m$ is vertex-transitive,
$\mathcal{G}_\alpha(x,y)$ depends only on the distance
$\min(|y-x|,m-|y-x|)$ between $x$ and $y$.  Therefore define
$\mathcal{G}_\alpha(a):=\mathcal{G}_\alpha(x,y)$, where $a=
|y-x|$; this induces the additional relation
$\mathcal{G}_\alpha(a) = \mathcal{G}_\alpha(m-a)$ for $1\leq a\leq
m-1$. From (\ref{eqn:normalizedGAlphaPseudoInverseRelations}), we have
\begin{eqnarray}
\chi(x=y)-\frac{1}{m} & = & \left(\mathcal{L}_S +
    \alpha\right)\mathcal{G}_\alpha(x,y) \nonumber \\
& = & \frac{1}{2}\left(2(1+\alpha)\mathcal{G}_\alpha(x,y)
    -\mathcal{G}_\alpha(x+1,y) - \mathcal{G}_\alpha(x-1,y)\right)
    \nonumber \\
& = & \frac{1}{2}\left((r+r^{-1})\mathcal{G}_\alpha(x,y)
    -\mathcal{G}_\alpha(x+1,y)-\mathcal{G}_\alpha(x-1,y)\right),
    \nonumber
\end{eqnarray}
where $\mathcal{L}_S$ is the normalized Laplacian of $C_m$.  We
can rewrite this as
$$ \mathcal{G}_{\alpha}(x-1,y)-r\mathcal{G}_{\alpha}(x,y)
    \ = \ \frac{2}{m} - 2\chi(x=y) + \frac{1}{r}(
    \mathcal{G}_{\alpha}(x,y)-r\mathcal{G}_{\alpha}(x+1,y)),
$$
which for $a>0$ becomes
\begin{eqnarray}
\mathcal{G}_\alpha(a+1)-r\, \mathcal{G}_\alpha(a) & = &
    \frac{2}{m} + \frac{1}{r}\left(\mathcal{G}_\alpha(a)
    -r\, \mathcal{G}_\alpha(a-1)\right) \nonumber \\
& = & \qquad \vdots \nonumber \\
& = & \frac{2}{m}
    + \cdots
    + \frac{1}{r^{a-1}}\frac{2}{m} + \frac{1}{r^a}
    \left(\mathcal{G}_\alpha(1)-r\, \mathcal{G}_\alpha(0)\right).
    \label{eqn:GAlphaCycleADiff}
\end{eqnarray}
When $a=0$, 
$$ \mathcal{G}_\alpha(1) - r\, \mathcal{G}_\alpha(0) =
    \frac{2}{m} - 2 + \frac{1}{r}\, (\mathcal{G}_\alpha(0)
    -r\, \mathcal{G}_\alpha(m-1)),$$
and since $\mathcal{G}_\alpha(1)=\mathcal{G}_\alpha(m-1)$, we
obtain
\begin{eqnarray}
\mathcal{G}_\alpha(1) & = & \frac{1}{m} - 1 + \frac{r+r^{-1}}{2}\,
    \mathcal{G}_\alpha(0) \nonumber \\
\mathcal{G}_{\alpha}(1)-r\mathcal{G}_{\alpha}(0) & = &
    \frac{1}{m} - 1 +
    \frac{-r+r^{-1}}{2}\, \mathcal{G}_\alpha(0)
    \label{eqn:GAlphaCycle1From0}
\end{eqnarray}
From this point, the reader who wishes to verify details is
encouraged to employ any standard computer algebra system.
Equations (\ref{eqn:GAlphaCycleADiff}) and (\ref{eqn:GAlphaCycle1From0})
define a recurrence with an initial condition on differences of
$\mathcal{G}_{\alpha}$, which is resolved by substituting
$\mathcal{G}_{\alpha}(1)-r\mathcal{G}_{\alpha}(0)$ from
(\ref{eqn:GAlphaCycle1From0}) into (\ref{eqn:GAlphaCycleADiff})
 and simplifying the geometric series.   For $a\geq 0$ this
yields
\begin{equation}\label{eqn:GAlphaCycleDAlpha}
\mathcal{G}_\alpha(a+1)-r\, \mathcal{G}_\alpha(a) =
    \frac{2}{r^{a-1}m} \frac{1-r^a}{1-r} + \frac{1}{r^a}
    \left(\frac{1}{m} - 1 +
    \frac{-r+r^{-1}}{2}\, \mathcal{G}_\alpha(0)\right).
\end{equation}
Denoting the right-hand side of (\ref{eqn:GAlphaCycleDAlpha}) by
$D_\alpha(a)$, for $a>0$ we have
\begin{eqnarray}
\mathcal{G}_\alpha(a) & = & r\, \mathcal{G}_\alpha(a-1) +
    D_\alpha(a-1) \nonumber \\
& = & \quad \vdots \nonumber \\
& = & r^a\mathcal{G}_\alpha(0) + r^{a-1} D_\alpha(0)
    + r^{a-2} D_\alpha(1) + \cdots + r^0 D_\alpha(a-1).
    \label{eqn:GAlphaCycleGFromD}
\end{eqnarray}
A careful but straightforward summing of geometric series in
(\ref{eqn:GAlphaCycleGFromD}) yields, for $a>0$,
\begin{eqnarray}
\mathcal{G}_\alpha(a) & = & \frac{1}{2}\, \mathcal{G}_\alpha(0)
    \frac{1+r^{2a}}{r^a} \nonumber \\
& & + \frac{2}{m} \frac{1}{r^{a-2}}
    \frac{1-r^a}{1-r} \left( \frac{1+r^a}{1-r^2} +
    \frac{r^{a-1}}{1-r} + \frac{1+r^a}{1+r} \frac{1-m}{2r}
    \right). \label{eqn:GAlphaCycleGFrom0}
\end{eqnarray}
Now using (\ref{eqn:GAlphaCycleGFrom0}), we set
$\mathcal{G}_\alpha(1)=\mathcal{G}_\alpha(m-1)$ and solve for
$\mathcal{G}_\alpha(0)$, obtaining
\begin{equation}\label{eqn:GAlphaCycleG0}
\mathcal{G}_\alpha(0) = -\frac{2}{m} \frac{r}{(r-1)^2} +
    \frac{2r(1+r^m)}{(r^2-1)(r^m-1)},
\end{equation}
which together with (\ref{eqn:GAlphaCycleGFrom0}) and
simplification yields
\begin{eqnarray*}
\mathcal{G}_\alpha(a) & = & \mathcal{G}_\alpha(0) -
    \frac{2\left(r^{a/2}-r^{-a/2}\right)\left(r^{m/2-a/2}-r^{-(m/2-a/2)}\right)}
    {\left(r-r^{-1}\right)\left(r^{m/2}-r^{-m/2}\right)} \nonumber\\
& = & \frac{-2}{m\left(r+r^{-1}-2\right)} +
    \frac{2\left(r^{m/2}+r^{-m/2}\right)}
    {\left(r-r^{-1}\right)\left(r^{m/2}-r^{-m/2}\right)} \nonumber
    \\
& & -  \frac{2\left(r^{a/2}-r^{-a/2}\right)
    \left(r^{m/2-a/2}-r^{-(m/2-a/2)}\right)}
    {\left(r-r^{-1}\right)\left(r^{m/2}-r^{-m/2}\right)}
    \nonumber \\
& = & \frac{-2}{m\left(r+r^{-1}-2\right)} +
    \frac{2(r^{m/2-|x-y|} +
    r^{-m/2+|x-y|})} {(r-r^{-1})(r^{m/2}-r^{-m/2})}\,
    . \nonumber
\end{eqnarray*}
Substituting $|y-x|$ for $a$ gives the desired formula for
$\mathcal{G}_{\alpha}$.
\end{proof}

By the definition of $\alpha$ and $r$, we may use the substitution
$r=e^{i\theta}$ to rewrite $(r^z+r^{-z})/2=\cos{z\theta}$ and
$(r^z-r^{-z})/2i=\sin{z\theta}$.  Together with the definition of
the Chebyshev polynomials of the first and second kinds, $T_n$ and
$U_n$, respectively; i.e.,
\begin{align*}
T_n(x) & := \cos{n\theta}\quad\mbox{and} \\
U_n(x) & := \frac{\sin{(n+1)\theta}}{\sin{\theta}},
\end{align*}
where $x=\cos{\theta}$, we obtain the following corollary to
Theorem \ref{thm:normalizedGAlphaCycleR}.

\begin{cor} \label{cor:normalizedGAlphaCycle}
Let $m\geq 3$.  For complex $\alpha\neq 0$ and  \ 
$0\leq x, y\leq m-1$, the generalized
Green's function $\mathcal{G}_\alpha$ for the cycle $C_m$ with
vertices \ $0,1,\ldots,m-1$ satisfies
\begin{eqnarray*}
\mathcal{G}_\alpha(x,y) & = & -\frac{1}{m\alpha}
    +\frac{T_{m/2-|y-x|}(1+\alpha)}{\alpha(2+\alpha)U_{m/2-1}(1+\alpha)}
    \nonumber
\end{eqnarray*}
where $T$ and $U$ are the Chebyshev polynomials of the first and
second kinds, respectively.
\end{cor}

\begin{proof}[Proof of Theorem \ref{thm:normalizedGTorus}:]
The theorem follows by applying the integral formula for products
of graphs without boundary in Corollary
\ref{cor:normalizedGProductWOBoundary} to the torus, where
$\Gamma=C_m$, $\Gamma'=C_n$, the $\phi'$'s are the orthonormal
basis described in Lemma \ref{lem:cycleMixedEigensystem}, 
$\mathcal{G}$
and $\mathcal{G}'$ are given by Theorem
\ref{thm:normalizedGCycle}, and $\mathcal{G}_{\alpha}$ is given by
Corollary \ref{cor:normalizedGAlphaCycle}.
\end{proof}

Combining Theorem \ref{thm:normalizedGTorus} with
(\ref{eqn:normalizedGFourier}) 
using the orthonormal eigensystem of 
Lemma \ref{lem:cycleMixedEigensystem} for both $C_m$ and $C_n$ yields the
following nontrivial identity.

\begin{cor}\label{cor:normalizedGTorusIdentity}
Let $m,n\geq 3$; $0\leq x, y\leq m-1$; and \ $0\leq x',y'\leq n-1$.  Then
\begin{align*}
\frac{1}{mn}\sum_{(j,k)\neq(0,0)}&
    \frac{\exp\big((2\pi ij/m)(y-x)\big)
    \exp\big((2\pi ik/n)(y'-x')\big)}
    {\big(1-\cos{\left(2\pi j/m\right)}/2-\cos{\left(2\pi
    k/n\right)}/2\big)}
        \nonumber \\
&= \ \frac{2}{n}\sum_{k=1}^{n-1}
    \exp\big((2\pi ik/n)(y'-x')\big)
    \left[-\frac{1}{m(1-\cos{(2\pi k/n)})}\right. \nonumber \\
& \left.+\frac{T_{m/2-|y-x|}(2-\cos{(2\pi k/n)})}
    {(1-\cos{(2\pi k/n)})
    (3-\cos{(2\pi k/n)})U_{m/2-1}(2-\cos{(2\pi k/n)})}
    \right]\nonumber \\
& +\frac{(m+1)(m-1)}{3mn} - \frac{|y-x|}{n} + \frac{(y-x)^2}{mn}
    \nonumber \\
& +\frac{(n+1)(n-1)}{3mn} - \frac{|y'-x'|}{m} +
    \frac{(y'-x')^2}{mn}, \nonumber
\end{align*}
where $T$ and $U$ are the Chebyshev polynomials of the first and
second kinds, respectively.
\end{cor}

The Laplacian of $C_m$ has a 1-dimensional eigenspace
corresponding to eigenvalue 0, and a second 1-dimensional
eigenspace corresponding to eigenvalue 2 iff $C_m$ is bipartite
(when $m$ is even).  Otherwise, all eigenspaces are 2-dimensional,
since $\lambda_j=\lambda_{m-j}$ for all $1\leq j\leq m-1$.  This
means that Corollary \ref{cor:normalizedGTorusIdentity} is only
one of a class of identities constructed by choosing orthonormal
eigensystems for $C_m$ and $C_n$ for the left-hand side, and a
possibly distinct orthonormal eigensystem for $C_n$ on the
right-hand side.

\subsection{The $t$-dimensional torus, 
$C_{m_1}\times C_{m_2}\times\cdots\times C_{m_t}$}

The bottleneck in computing the normalized Green's function for 
the $t$-torus via Corollary \ref{cor:normalizedGProductWOBoundary} 
is the lack of a formula for the generalized Green's function
for any $t\geq 2$.  Thus the decomposition into a cycle,
for which $\mathcal{G}_{\alpha}$ is given by Corollary 
\ref{cor:normalizedGAlphaCycle}, and a 
$(t-1)$-torus is required.

Before giving the normalized Green's function of $C_{m_1}\times
\cdots\times C_{m_t}$, we present the information on the
components still needed.  Choose and label the eigensystem of each
$C_{m_s}$ by
$$\{(\lambda^{(s)}_{j_s},\phi^{(s)}_{j_s}):0\leq j_s\leq m_s-1\},$$
where $m_s\geq 3$, $0=\lambda^{(s)}_{0}<\lambda^{(s)}_{1}\leq
\cdots \leq \lambda^{(s)}_{m_s-1}$, and the vectors
$\{\phi^{(s)}_{j_s}:0\leq j_s\leq m_s-1\}$ are orthonormal.  The
eigenvalues of $C_{m_2}\times \cdots\times C_{m_t}$ are averages
of the eigenvalues of the factors $C_{m_s}$, and the corresponding
eigenvectors are products of the eigenvectors of the factors. This
is summarized in the next well-known lemma, whose proof is a straightforward
induction on (\ref{eqn:productEigensystem}).
\begin{lemma}\label{lem:normalizedGRepeatedCycleProductEigensystem}
The eigenvalues of\, $C_{m_2}\times \cdots\times C_{m_t}$ are
$$  \Lambda_{j_2,\ldots,j_t} \ = \
\frac{\lambda^{(2)}_{j_2}+\cdots+\lambda^{(t)}_{j_t}}{t-1},$$
where $0\leq j_s \leq m_s-1$ for all $2\leq s\leq t$, with
corresponding eigenvectors
$$\Phi_{j_2,\ldots,j_t}(x_2,\ldots,x_t) \ := \
\prod_{s=2}^{t}\phi^{(s)}_{j_s}(x_s).$$
\end{lemma}

For the following theorem, let $\mathcal{G}$ be the normalized
Green's function for $C_{m_1}$ from Theorem
\ref{thm:normalizedGCycle}, and $\mathcal{G}'$ the normalized
Green's function for $C_{m_2}\times\cdots C_{m_{t}}$.

\begin{theorem}\label{thm:normalizedGRepeatedCycleProduct}
Let\, $t\geq 2$. Let\, $0\leq x_{j_s},y_{j_s}\leq m_s-1$ where
$m_s\geq 3$ for $1\leq s\leq t$. The $t$-dimensional torus
$C_{m_1}\times\cdots\times C_{m_t}$ has normalized Green's
function
\begin{align}
\mathbf{G}((&x_1,\ldots,x_t),(y_1,\ldots,y_t)) \ = \
    t\sum_{K\neq(0,\ldots,0)}\Phi_K(x_2,\ldots,x_{t})
    \overline{\Phi_K(y_2,\ldots,y_{t}})
    \mathcal{G}_{\Lambda_K}(x_1,y_1) \nonumber \\
& \quad +\frac{t}{(t-1)m_1}\,
    \mathbf{G}'((x_2,\ldots,x_{t}),(y_2,\ldots,y_{t}))
    +\frac{t}{m_2\cdots m_{t}}\, \mathcal{G}(x_1,y_1),
    \label{eqn:normalizedGTTorus}
\end{align}
where $K$ ranges over all indices $(j_2,\ldots,j_{t})\neq
(0,\ldots,0)$, $\Phi$ and $\Lambda$ are defined in Lemma
\ref{lem:normalizedGRepeatedCycleProductEigensystem},
$\mathcal{G}_\alpha$ is given by Corollary
\ref{cor:normalizedGAlphaCycle}, and $\mathbf{G}'$ is the 
normalized Green's function for $C_{m_2}\times\cdots\times C_{m_t}$.
\end{theorem}
\begin{proof} The proof proceeds by using $\Gamma=C_{m_1}$ and
$\Gamma'=C_{m_2}\times\cdots \times C_{m_{t}}$ in Corollary
\ref{cor:normalizedGProductWOBoundaryDDPrime}.  The degree of
$\Gamma$ is $d=2$, and the degree of $\Gamma'$ is $d'=2(t-1)$. The
result follows. \end{proof}

Although $\mathbf{G}'$ in Theorem
\ref{thm:normalizedGRepeatedCycleProduct} may not already be
known, it can be computed inductively from repeated applications
of the theorem.  Determination of $\mathcal{G}_\alpha$ for any
small product $C_1\times\cdots\times C_{t'}$ would allow the
reduction in the number of applications required by a factor of
$t'$. The following formula for the 3-torus is obtained by two
applications of Theorem \ref{thm:normalizedGRepeatedCycleProduct},
first taking the product of $\Gamma=C_m$ with $\Gamma'=C_m\times C_m$, and
then the product of $\Gamma=C_m$ with $\Gamma'=C_m$.

\begin{cor} For\, $0\leq x_1,y_1,x_2,y_2,x_3,y_3\leq m-1$ where $m\geq 3$,
the 3-dimensional torus\, $C_m\times C_m\times C_m$ has normalized
Green's function
\begin{align*}
\mathcal{G}^{\times}((x_1,&x_2,x_3),(y_1,y_2,y_3))
    \ = \ \frac{3}{m^2}
    \sum_{(j,k)\neq(0,0)}
    \Big[\exp\big((2\pi ij/m)(y_2-x_2)\big)
    \nonumber \\
& \quad \ \exp\big((2\pi ik/m)(y_3-x_3)\big)
    \mathcal{G}_{\left(1-\cos{(2\pi j/m)}/2-
    \cos{(2\pi k/m)}/2\right)}(x_1,y_1)\Big] \nonumber\\
& + \frac{3}{m^2}\sum_{j=1}^{m}
    \exp\big((2\pi ik/m)(y_3-x_3)\big)
    \mathcal{G}_{\left(1-\cos{(2\pi j/m)}\right)}(x_2,y_2) \nonumber\\
& +\frac{3}{m^2}\left(\frac{(m+1)(m-1)}{6m} - |y_3-x_3|
    +\frac{(y_3-x_3)^2}{m}\right) \nonumber \\
& +\frac{3}{m^2}\left(\frac{(m+1)(m-1)}{6m} - |y_2-x_2| +
    \frac{(y_2-x_2)^2}{m}\right)  \nonumber \\
& +\frac{3}{m^2}\left(\frac{(m+1)(m-1)}{6m} - |y_1-x_1| +
    \frac{(y_1-x_1)^2}{m}\right). \nonumber
\end{align*}
where $\mathcal{G}_\alpha$ is given in Corollary
\ref{cor:normalizedGAlphaCycle}.
\end{cor}

These compact formulas for Green's functions of tori offer
fast alternatives to computing pseudo-inverses of their Laplacians
directly.  This increase in speed, essentially due to the symmetry
of the torus, is reflected in $O(\log n)$ computational 
complexity of the Chebyshev polynomials $T_n$
and $U_n$. Various algorithms for computing $T_n$ and $U_n$
are given in \cite{K95}, and a more theoretical treatment of types of
polynomials computable in $O(\log n)$ appears in \cite{F89}. The
following corollary to Theorem
\ref{thm:normalizedGRepeatedCycleProduct} is significant because
the Laplacian ($\mathcal{L}$, $L$, or $\Delta$) of the torus on
$n$ vertices has rank $n-1$, and so computing its pseudo-inverse
provides along the way the inverse of an $(n-1)\times (n-1)$
matrix.

\begin{cor}
Matrix pseudo-inversion of the Laplacian of the $t$-dimensional
torus with $n$ vertices via its Green's function is \ $O(t\cdot
n^{2-1/t}\log n)$, provided that the matrix itself is not
completely reconstructed.
\end{cor}

\begin{proof}
We assume the $t$-dimensional torus is $C_{m_1}\times\cdots\times
C_{m_t}$, where $m_s\geq 3$ for $1\leq s\leq t$ and $\prod_{s=1}^t
m_s = n$. It suffices to compute one row of the pseudo-inverse of
the normalized Laplacian in order to know the entire inverse, due
to the translational symmetry of the torus; i.e., since
\begin{equation*}
\mathbf{G}((x_1,\ldots, x_t),(y_1,\ldots, y_t)) \ = \
\mathbf{G}((0,\ldots, 0),(|y_1-x_1|,\ldots, |y_t-x_t|)).
\end{equation*}
(In fact, because $|y_s-x_s|$ can be replaced by $m-|y_s-x_s| \
(\mathrm{mod} \ m)$ without changing the value of
$\mathbf{G}$, only $\prod_{s=1}^t\lceil m_s/2\rceil$ of
these entries must actually be computed.) Without loss of generality,
$m_1\geq \cdots \geq m_t$. Then the summation
term on the right-hand side of (\ref{eqn:normalizedGTTorus}) has
at most $n^{1-1/t}$ summands, which can each be computed in $O(t
\log{n})$ time. The second and third terms can also be computed in
$O(t\,n^{1-1/t}\log{n})$ time, and we must compute $n$ terms total
to know all entries of the pseudo-inverse of $\mathcal{L}$.
\end{proof}

For example, the time complexity of computing the Green's function
for the 3-torus with $n$ vertices is $O(n^{5/3}\log n)$.
Such a quick
pseudo-inversion formula is surprising, since matrix inversion in
general has the same complexity as matrix multiplication (see
\cite{BB88}). Matrix multiplication is known to be
$O(n^\omega)$, where $2\leq \omega\leq 2.376$ (see \cite{CW90}).
Surprisingly, for large $n$ we can compute all of the values for
the pseudo-inverse of the Laplacian for the torus faster than we
can write down the whole matrix. Of course, requiring the
presentation of the whole matrix rather than just the first row
increases the complexity to $O(t\,n^{1-1/t}\log{n}+n^2)$.

\section{Hitting times from Green's functions\label{sec:hittingTime}}

The equivalence of Green's functions to the fundamental matrix $Z$
of (\ref{eqn:fundMatFormal}) under similarity transformation
allows many quantities in random walks to be computed using $\mathcal{G}$.
The {\em hitting time} $Q(x,y)$ of a simple random walk starting at
vertex $x$ with target vertex $y$ is the expected number of steps
to reach vertex $y$ for the first time by starting at $x$ and at each step
moving to any neighbor of $x$ with equal probability. 
In \cite{CY00}, Chung and Yao show the
following relationship between $\mathcal{G}$ and $Q$.

\begin{theorem}[Chung, Yao]
The hitting time $Q(x,y)$ satisfies
$$  Q(x,y) = \frac{vol(\Gamma)}{d_y}\mathcal{G}(y,y)
    -\frac{vol(\Gamma)}{\sqrt{d_xd_y}}\mathcal{G}(x,y). $$
\end{theorem}

\begin{figure}[h]
\centerline{\includegraphics[width=251pt,height=131pt]{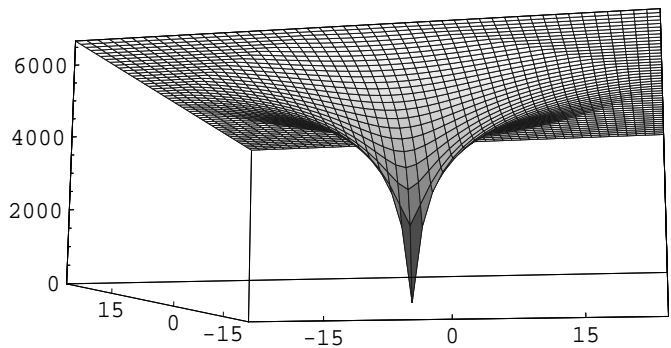}} 
\caption{Graph of the hitting time of a random walk on
the torus $C_{49}\times C_{49}$, laid out as a square with periodic
boundary} \label{fig:torusHitting}
\end{figure}

We immediately have a computational formula for the hitting time
whenever $\mathcal{G}$ is known.  For instance, whenever $\Gamma$
is regular with $n$ vertices,
\begin{equation}\label{eqn:hittingGreen}
Q(x,y) = n\left(\mathcal{G}(y,y)
    -\mathcal{G}(x,y)\right).
\end{equation}
Figure \ref{fig:torusHitting} plots hitting time in the case of 
$C_{49}\times C_{49}$, using Theorem \ref{thm:normalizedGTorus} and 
(\ref{eqn:hittingGreen}).
The domain consists of the vertices of the torus 
laid out in a square grid with periodic boundary,
making opposite ends
adjacent.  The vertical axis plots the hitting time of a random
walk initiated at $(0,0)$ with target $(x,y)$, achieving a 
minimum of 0 at $(0,0)$ and leveling off at just above
6000 steps to reach 
vertices farthest from the start of the walk.
If (\ref{eqn:hittingGreen}) is computed using techniques such as
Corollary \ref{cor:normalizedGAlphaCycle}, the hitting time
expression will involve orthogonal polynomials.  Aldous and Fill
claim in \cite{AF02} that orthogonal polynomials appear
whenever the graph has sufficient symmetry, but the dependence
remains largely unstudied.

\section*{Acknowledgment}\vskip -.1in
The author would like to give special thanks to Fan R.\ K.\ Chung,
who provided assistance in the form of many valuable discussions
on the material in this paper.


\begin{thebibliography}{10}

\bibitem{AF02}
D.~Aldous and J.~Fill.
\newblock Reversible markov chains and random walks on graphs.
\newblock Manuscript, http://www.stat.berkeley.edu/users/aldous/RWG/book.html.


\bibitem{BB88}
Gilles Brassard and Paul Bratley.
\newblock {\em Algorithmics: Theory and Practice}.
\newblock Prentice Hall Inc., Englewood Cliffs, NJ, 1988.

\bibitem{CY00}
Fan Chung and S.-T. Yau.
\newblock Discrete {G}reen's functions.
\newblock {\em J. Combin. Theory Ser. A}, 91(1-2):191--214, 2000.

\bibitem{C97}
Fan R.~K. Chung.
\newblock {\em Spectral graph theory}, volume~92 of {\em CBMS Regional
  Conference Series in Mathematics}.
\newblock AMS Publications, 1997.

\bibitem{CW90}
Don Coppersmith and Shmuel Winograd.
\newblock Matrix multiplication via arithmetic progressions.
\newblock {\em J. Symbolic Comput.}, 9(3):251--280, 1990.

\bibitem{D79}
Philip~J. Davis.
\newblock {\em Circulant matrices}.
\newblock John Wiley \& Sons, New York-Chichester-Brisbane, 1979.

\bibitem{DS84}
Peter~G. Doyle and J.~Laurie Snell.
\newblock {\em Random walks and electric networks}, volume~22 of {\em Carus
  Mathematical Monographs}.
\newblock Mathematical Association of America, Washington, DC, 1984.

\bibitem{E02}
Robert~B. Ellis.
\newblock {\em Chip-firing games with Dirichlet eigenvalues and discrete
  Green's functions}.
\newblock PhD thesis, University of California at San Diego, 2002.

\bibitem{F89}
R.~Fateman.
\newblock Lookup tables, recurrences and complexity.
\newblock In {\em Proceedings of ISSAC 89}, pages 68--73, New York, 1989. ACM
  Press.

\bibitem{K95}
W.~Koepf.
\newblock Efficient computation of orthogonal polynomials in computer algebra.
\newblock Preprint SC 95-42, December 1995.

\bibitem{T00}
Andr{\'a}s Telcs.
\newblock Transition probability estimates for reversible {M}arkov chains.
\newblock {\em Electron. Comm. Probab.}, 5:29--37, 2000.

\end{thebibliography}
\end{document}